\documentclass[12pt]{amsart}
\usepackage{amsfonts,amssymb,amscd,amsmath,enumerate,verbatim}
\usepackage[latin1]{inputenc}
\usepackage{amscd}
\usepackage{latexsym}
\usepackage{mathptmx}
\usepackage{multicol}
\input xy
\xyoption{all}
\def\NZQ{\mathbb}               % the font for N,Z,Q,R,C

\def\ZZ{{\NZQ Z}}

\def\frk{\mathfrak}               % font for "Fraktur"

\def\Phi{{\frk N}}
\def\opn#1#2{\def#1{\operatorname{#2}}} % to make operators
\opn\chara{char} 
\opn\length{\ell} 
\opn\pd{pd} 
\opn\rk{rk}
\opn\projdim{proj\,dim} 
\opn\injdim{inj\,dim} 
\opn\rank{rank}
\opn\depth{depth} 
\opn\grade{grade} 
\opn\height{height}
\opn\embdim{emb\,dim} 
\opn\codim{codim}

\opn\Tr{Tr} 
\opn\bigrank{big\,rank}
\opn\superheight{superheight}
\opn\lcm{lcm}
\opn\trdeg{tr\,deg}%\emph{
\opn\reg{reg} 
\opn\lreg{lreg} 
\opn\ini{in} 
\opn\lpd{lpd}
\opn\size{size}
\opn\mult{mult}
\opn\dist{dist}
\opn\cone{cone}
\opn\lex{lex}
\opn\rev{rev}
\opn\im{im}
\opn\m{m}
\opn\div{div} \opn\Div{Div} \opn\cl{cl} \opn\Cl{Cl}
\opn\Spec{Spec} \opn\Supp{Supp} \opn\supp{supp} \opn\Sing{Sing}
\opn\Ass{Ass} \opn\Min{Min}
\opn\Ann{Ann} \opn\Rad{Rad} \opn\Soc{Soc}
\opn\Syz{Syz} \opn\Im{Im} \opn\Ker{Ker} \opn\Coker{Coker}
\opn\Am{Am} \opn\Hom{Hom} \opn\Tor{Tor} \opn\Ext{Ext}
\opn\End{End} \opn\Aut{Aut} \opn\id{id} \opn\ini{in}

\opn\nat{nat}
\opn\pff{pf}%   \pf exists already
\opn\Pf{Pf} \opn\GL{GL} \opn\SL{SL} \opn\mod{mod} \opn\ord{ord}
\opn\Gin{Gin}
\opn\Hilb{Hilb}\opn\adeg{adeg}\opn\std{std}\opn\ip{infpt}
\opn\Pol{Pol}
\opn\sat{sat}
\opn\Var{Var}
\opn\Gen{Gen}

\opn\aff{aff} \opn\con{conv} \opn\relint{relint} \opn\st{st}
\opn\lk{lk} \opn\cn{cn} \opn\core{core} \opn\vol{vol}
\opn\link{link} \opn\star{star}
\opn\gr{gr}

\def\pot#1#2{#1[\kern-0.28ex[#2]\kern-0.28ex]}

\opn\dirlim{\underrightarrow{\lim}}
\opn\inivlim{\underleftarrow{\lim}}
\let\to=\rightarrow

\def\Implies{\ifmmode\Longrightarrow \else
        \unskip${}\Longrightarrow{}$\ignorespaces\fi}
\def\implies{\ifmmode\Rightarrow \else
        \unskip${}\Rightarrow{}$\ignorespaces\fi}
\def\iff{\ifmmode\Longleftrightarrow \else
        \unskip${}\Longleftrightarrow{}$\ignorespaces\fi}

\let\:=\colon
\newtheorem{Theorem}{Theorem}[section]
\newtheorem{Lemma}[Theorem]{Lemma}
\newtheorem{Corollary}[Theorem]{Corollary}

\newtheorem{Remark}[Theorem]{Remark}

\newtheorem{Example}[Theorem]{Example}

\let\epsilon\varepsilon
\let\phi=\varphi
\let\kappa=\varkappa
\def\qed{\ifhmode\textqed\fi
      \ifmmode\ifinner\quad\qedsymbol\else\dispqed\fi\fi}
\def\textqed{\unskip\nobreak\penalty50
       \hskip2em\hbox{}\nobreak\hfil\qedsymbol
       \parfillskip=0pt \finalhyphendemerits=0}
\def\dispqed{\rlap{\qquad\qedsymbol}}

\opn\dis{dis}
\opn\height{height}
\opn\dist{dist}
\def\pnt{{\raise0.5mm\hbox{\large\bf.}}}

\opn\Lex{Lex}

\begin{document}
\title{Induced matching numbers of finite graphs and edge ideals}
\author{Takayuki Hibi, Hiroju Kanno and Kazunori Matsuda}
\address{Takayuki Hibi,
Department of Pure and Applied Mathematics,
Graduate School of Information Science and Technology,
Osaka University, Suita, Osaka 565-0871, Japan}
\email{hibi@math.sci.osaka-u.ac.jp}
\address{Hiroju Kanno, 
Department of Pure and Applied Mathematics,
Graduate School of Information Science and Technology,
Osaka University, Suita, Osaka 565-0871, Japan}
\email{u825139b@ecs.osaka-u.ac.jp}
\address{Kazunori Matsuda,
Department of Pure and Applied Mathematics,
Graduate School of Information Science and Technology,
Osaka University, Suita, Osaka 565-0871, Japan}
\email{kaz-matsuda@ist.osaka-u.ac.jp}
% \thanks{
% }
\subjclass[2010]{05C69, 05C70, 05E40, 13D40, 13H10}
\keywords{Castelnuovo--Mumford regularity, $h$-polynomial, edge ideal, induced matching number.}
\begin{abstract}
Let $G$ be a finite simple graph on the vertex set $V(G) = \{x_1, \ldots, x_n\}$ and $I(G) \subset K[V(G)]$ its edge ideal, where $K[V(G)]$ is the polynomial ring in $x_1, \ldots, x_n$ over a field $K$ with each $\deg x_i = 1$ and where $I(G)$ is generated by those squarefree quadratic monomials $x_ix_j$ for which $\{x_i, x_j\}$ is an edge of $G$.  In the present paper, given integers $1 \leq a \leq r$ and $s \geq 1$, the existence of a finite connected simple graph $G = G(a, r, d)$ with $\im(G) = a$, $\reg(R/I(G)) = r$ and $\deg h_{K[V(G)]/I(G)} (\lambda) = s$, where $\im(G)$ is the induced matching number of $G$ and where $h_{K[V(G)]/I(G)} (\lambda)$ is the $h$-polynomial of $K[V(G)]/I(G)$. 
\end{abstract}
\maketitle
\section*{Introduction}
The recent papers \cite{HM1}, \cite{HM2}, and \cite{HMVT} study the relation between regularity of monomial ideals and degree of their $h$-polynomials.  The present paper follows these previous work and develop the combinatorial aspect on regularity of edge ideals of finite simple graphs.
    
Let $R = K[x_1, \ldots, x_n]$ denote the polynomial ring in $n$ variables over a field $K$ with each $\deg x_i = 1$ and $I \subset R$ a homogeneous ideal of $R$ with $\dim R/I = d$.  The Hilbert series $H_{R/I}(\lambda)$ of $R/I$ is of the form 
\[
H_{R/I}(\lambda) = \left( h_{0}(R/I) + h_{1}(R/I)\lambda + h_{2}(R/I)\lambda^2 + \cdots + h_{s}(R/I)\lambda^s \right)/(1 - \lambda)^{d},
\] 
where each $h_{i}(R/I) \in \ZZ$ (\cite[Proposition 4.4.1]{BH}).  
We say that 
\[
h_{R/I}(\lambda) = h_{0}(R/I) + h_{1}(R/I)\lambda + h_{2}(R/I)\lambda^2 + \cdots + h_{s}(R/I)\lambda^s
\] 
with $h_{s}(R/I) \neq 0$ and $h_{R/I}(1) \neq 0$  is the {\em $h$-polynomial} of $R/I$.  
Let $\reg(R/I)$ denote the ({\em Castelnuovo--Mumford}\,) {\em regularity} \cite[p.~168]{BH} of $R/I$. 

Let $G = (V(G), E(G))$ be a finite simple graph on the vertex set $V(G) = \{x_{1}, \ldots, x_{n}\}$ with the edge set $E(G)$.  Let $R = K[V(G)] = K[x_{1}, \ldots, x_{n}]$.  The {\em edge ideal} of $G$ is the ideal 
\[
I(G) = \left( x_{i}x_{j} : \{x_{i}, x_{j}\} \in E(G) \right) \subset R. 
\] 

%A set of edges $\{e_{1}, \ldots, e_{s} \} \subset E(G)$ is said to be a {\em matching} of $G$ 
%if $e_{i} \cap e_{j} \neq \emptyset$ for all $1 \leq i < j \leq s$. 
A subset $M = \{e_{1}, \ldots, e_{s} \} \subset E(G)$ is said to be a {\em matching} of $G$ if, for all $e_{i}$ and $e_{j}$ with $i \neq j$ belonging to $M$, one has $e_i \cap e_j = \emptyset$.  The {\em matching number} $\m(G)$ of $G$ is the maximum cardinality of the matchings of $G$.   
A matching $M = \{e_{1}, \ldots, e_{s} \} \subset E(G)$ is said to be an {\em induced matching} of $G$ if for all $e_{i}$ and $e_{j}$ with $i \neq j$ belonging to $M$, there is no edge $e \in E(G)$ with $e \cap e_{i} \neq \emptyset$ and $e \cap e_{j} \neq \emptyset$.   
%The {\em matching number} $m(G)$ of $G$ is the maximum cardinality 
%of the matchings of $G$. 
The {\em induced matching number} $\im(G)$ of $G$ is the maximum cardinality 
of the induced matchings of $G$.  It is known (\cite[Theorem 6.7]{HaVanTuyl} and \cite[Lemma 2.2]{K}) that 
\[
\im(G) \leq \reg \left(R/I(G)\right) \leq \m(G).
\] 
We refer the reader to \cite{BC1}, \cite{HHKT}, \cite{Ku} and \cite{W} for further information about (induced) matching numbers and regularity of edge ideals.

In the current trend on the study of regularity of powers of edge ideals, gap free graphs play an important role (\cite{BC2}, \cite{CKV} and \cite{DHS}).  Recall that a {\em gap free} graph is a finite simple graph $G$ such that, for edges $e$ and $e'$ of $G$ with $e \cap e' = \emptyset$, there is an edge $e''$ of $G$ with $e \cap e'' \neq \emptyset$ and $e' \cap e'' \neq \emptyset$.  In other words, a finite simple graph $G$ is gap free if and only if $\im(G) = 1$.

Now, the final goal of the present paper is to show that 
  
\begin{Theorem}
\label{main}
Given integers $1 \leq a \leq r$ and $s \geq 1$, there exists a finite connected simple graph $G = G(a, r, s)$ with
\[
\im(G) = a, \, \, \, \, \, \reg(R/I(G)) = r, \, \, \, \, \, \deg h_{R/I(G)} (\lambda) = s. 
\]
\end{Theorem} 

Theorem \ref{main} can be an extensive generalization of \cite[Theorem 3.1]{HMVT}.  As one of the direct consequences of Theorem \ref{main}, it follows that

\begin{Corollary}
\label{gapfree}
Given integers $r \geq 1$ and $s \geq 1$, there exists a gap free graph $G = G(r, s)$ with $\reg(R/I(G)) = r$ and $\deg h_{R/I(G)} (\lambda) = s$. 
\end{Corollary}

%\section{Proof of Theorem \ref{main}}
\section{Preparation for Theorem \ref{main}}

In order to prove our Theorem \ref{main}, we will prepare several lemmata.
Let, as before, $R = K[x_1, \ldots, x_n]$ denote the polynomial ring in $n$ variables 
over a field $K$ with each $\deg x_i = 1$.  
From the definition of regularity in terms of graded Betti numbers (\cite[p.~48]{HH}), Lemma \ref{LemA} below follows immediately.  

\begin{Lemma}\label{LemA}
Let $I \subset R$ be a proper homogeneous ideal.  Then $\reg (R/I) = \reg (I) - 1$.  
\end{Lemma}

\begin{Lemma}[{\cite[Lemma 3.2]{HT}}]
\label{LemB}%$($\cite[Lemma 3.2]{HT}$)$
Let $R_{1} = K[x_{1}, \ldots, x_{m}]$ and $R_{2} = K[y_{1}, \ldots, y_{n}]$ be polynomial rings over a field $K$. 
Let $I_{1}$ be a nonzero homogeneous ideal of $R_{1}$ and $I_2$ that of $R_2$.  
Let $R = R_{1} \otimes_{K} R_{2} = K[x_{1}, \ldots, x_{m}, y_{1}, \ldots, y_{n}]$ and regard $I_{1} + I_{2}$ as a homogeneous ideal of $R$. Then
\begin{enumerate}%[$(1)$]
	%\item[$(1)$] $\reg (I_{1}I_{2}) = \reg (I_{1}) + \reg (I_{2})$\,$;$ 
	\item[$(1)$] $\reg (I_{1} + I_{2}) = \reg (I_{1}) + \reg (I_{2}) - 1$. 
	\item[$(2)$] $\reg (R/(I_{1} + I_{2})) = \reg (R_{1}/I_{1}) + \reg (R_{2}/I_{2})$. 
	%\item[$(4)$] $H_{R/(I_{1} + I_{2})}(\lambda) = H_{R_{1}/(I_{1})}(\lambda) + H_{R_{2}/(I_{2})}(\lambda)$.    
\end{enumerate}
\end{Lemma}

\begin{Lemma}[{\cite[Lemma 2.10]{DHS}}]
\label{LemC}
Let $I \subset R$ be a squarefree monomial ideal and $x_i$ a variable of $R$ which appears in 
a monomial belonging to the unique minimal system of monomial generators of $I$.  
Then $\reg (I) = \reg (I : (x)) + 1$ or  $\reg (I) = \reg (I + (x))$. 
\end{Lemma}

Let us recall the definition of an independent set. 
Let $G = (V(G), E(G))$ be a finite simple graph on the vertex set 
$V(G) = \{x_{1}, \ldots, x_{n}\}$ with the edge set $E(G)$. 
Given any subset $W \subset V(G)$, the {\em induced subgraph} of $G$ 
on $W$ is the graph $G_{W} = (W, E(G_{W}))$, 
where $E(G_{W}) = \{ e = \{x_{i}, x_{j}\} : i, j \in W \}$. 
A set of vertices $W \subset V(G)$ is an {\em independent set} if 
$\{ x_{i}, x_{j} \} \not\in E(G)$ for all $x_{i}, x_{j} \in W$. 
In particular, the empty set $\emptyset$ is an independent set for all finite graph $G$. 

Let $G$ be a finite simple graph on $V(G) = \{ x_{1}, \ldots, x_{n} \}$. 
Suppose that $G$ has no isolated vertex. 
Let $S \subset V(G)$ be an independent set of $G$. 
Note that $0 \leq |S| \leq d = \dim R/I(G)$ since $\dim R/I(G)$ is equal to the maximum cardinality 
of independent sets of $G$. 
The graph $G^{S}$ is defined by
\begin{itemize}
	\item $V(G^{S}) = V(G) \cup \{ x_{n + 1} \}$, where $x_{n + 1}$ is a new vertex. 
	\item $E(G^{S}) = E(G) \cup \left\{ \{ x_{i}, x_{n + 1} \} : x_{i} \not\in S \right\}$. 
\end{itemize}
We call $G^{S}$ the $S$-{\em suspension} of $G$. 
Note that $G^{\emptyset}$ coincides with the suspension \cite[p.141]{HNOS} of $G$ in usual sense. 

In addition, let $\{ x_{i}, x_{j} \} \in E(G)$ be an edge of $G$ and let $S \subset V(G)$ be an independent set of $G$ 
such that $\{x, x_{i}\} \not\in E(G)$ and $\{x, x_{j}\} \not\in E(G)$ for all $x \in S$. 
Then the graph $G^{\{x_{i}, x_{j}\}, S}$ is defined by
\begin{itemize}
	\item $V\left(G^{\{x_{i}, x_{j}\}, S}\right) = V(G) \cup \{ x_{n + 1} \}$, where $x_{n + 1}$ is a new vertex. 
	\item $E\left(G^{\{x_{i}, x_{j}\}, S}\right) = E(G) \cup \left\{ \{ x_{k}, x_{n + 1} \} : x_{k} \neq x_{i}, x_{k} \neq x_{j} \ \text{and} \ x_{k} \not\in S \right\}$. 
\end{itemize}
We call $G^{\{x_{i}, x_{j}\}, S}$ the $\{x_{i}, x_{j}\}, S$-{\em suspension} of $G$. 

\begin{Example}\normalfont
\label{C5}
Let $G = C_{5}$ be the pentagon graph with $V(G) = \{x_{1}, x_{2}, x_{3}, x_{4}, x_{5}\}$ and 
$E(G) = \left\{ \{x_{1}, x_{2}\}, \{x_{2}, x_{3}\}, \{x_{3}, x_{4}\}, \{x_{4}, x_{5}\}, \{x_{1}, x_{5}\} \right\}$. 
Then the graphs $G^{\emptyset}$ and $G^{\{x_{1}, x_{2}\}, \emptyset}$ are as follows. 

\begin{figure}[htbp]
  \centering

\bigskip

\begin{xy}
	\ar@{} (0,0);(50, 10) *++!D{x_{1}} *\dir<4pt>{*} = "A1";
	\ar@{-} "A1";(34, 0) *++!R{x_{5}} *\dir<4pt>{*} = "E1";
	\ar@{-} "A1";(66, 0) *++!L{x_{2}} *\dir<4pt>{*} = "B1";
	\ar@{-} "E1";(40, -18) *++!R{x_{4}} *\dir<4pt>{*} = "D1";
	\ar@{-} "D1";(60, -18) *++!L{x_{3}} *\dir<4pt>{*} = "C1";
	\ar@{-} "B1";"C1";
	\ar@{-} "A1";(50, -4) *++!U{x_{6}} *\dir<4pt>{*} = "F1";
	\ar@{-} "B1";"F1";
	\ar@{-} "C1";"F1";
	\ar@{-} "D1";"F1";
	\ar@{-} "E1";"F1";
	\ar@{} (0,0);(110, 10) *++!D{x_{1}} *\dir<4pt>{*} = "A2";
	\ar@{-} "A2";(94, 0) *++!R{x_{5}} *\dir<4pt>{*} = "E2";
	\ar@{-} "A2";(126, 0) *++!L{x_{2}} *\dir<4pt>{*} = "B2";
	\ar@{-} "E2";(100, -18) *++!R{x_{4}} *\dir<4pt>{*} = "D2";
	\ar@{-} "D2";(120, -18) *++!L{x_{3}} *\dir<4pt>{*} = "C2";
	\ar@{-} "B2";"C2";
	\ar@{} "A2";(110, -4) *++!U{x_{6}} *\dir<4pt>{*} = "F2";
	\ar@{} "B2";"F2";
	\ar@{-} "C2";"F2";
	\ar@{-} "D2";"F2";
	\ar@{-} "E2";"F2";
\end{xy}
  \caption{$G^{\emptyset}$ (left) and $G^{\{x_{1}, x_{2}\}, \emptyset}$ (right)}
%  \caption{$G_{V(G) \setminus \{v_{1}\}}$ (left) and $G_{V(G) \setminus N_{G}[v_{1}]}$ (right) in the case $m=n=1$ and $t_1 >0$}
 
\end{figure}
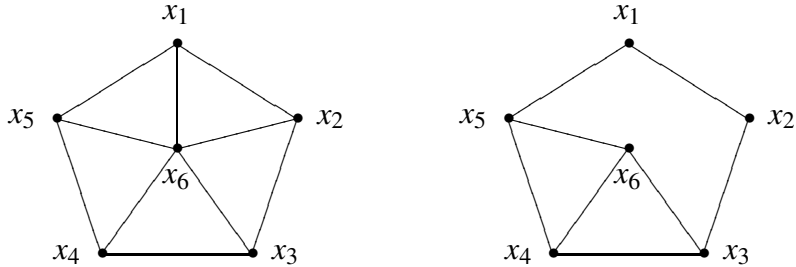

\bigskip

\end{Example}

Lemma \ref{S-suspension} below mentions the induced matching number, regularity and Hilbert series of 
the $S$-suspension of a graph. \\

\begin{Lemma}
\label{S-suspension}
Let $G$ be a finite simple graph on the vertex set $V(G) = \{x_{1}, \ldots, x_{n}\}$
such that $G$ has no isolated vertices. 
Let $S \subset V(G)$ be an independent set of $G$. 
Write $R' = R \otimes_{K} K[x_{n + 1}]$,  
and $h_{R/I(G)}(\lambda) =  h_{0}(R/I(G)) + h_{1}(R/I(G))\lambda + h_{2}(R/I(G))\lambda^2 + \cdots + h_{s}(R/I(G))\lambda^s$
with $h_{s}(R/I(G)) \neq 0$. 
Then one has 
\begin{enumerate}
	\item[$(1)$]  $\im(G^{S}) = \im(G)$. 
	\item[$(2)$]  $\reg(R'/I(G^{S})) = \reg(R/I(G))$. 
	\item[$(3)$]  The Hilbert series of $R'/I(G^{S})$ is 
	\begin{eqnarray*} H_{R'/I(G^{S})}(\lambda) &=& H_{R/I(G)}(\lambda) + \frac{\lambda}{(1 - \lambda)^{|S| + 1}} \\
	&=& \begin{cases} 
		\displaystyle \frac{h_{R/I(G)}(\lambda) + \lambda(1 - \lambda)^{d - |S| - 1} }{(1 - \lambda)^{d}} & \text{$(|S| \leq d - 1)$}, \\
		\displaystyle \frac{(1 - \lambda)h_{R/I(G)}(\lambda) + \lambda}{(1 - \lambda)^{d + 1}} & \text{$(|S| = d)$}. 
		\end{cases} 
	\end{eqnarray*}
	In particular, if $|S| = d - s$, we have 
	\begin{enumerate}
		\item[$(i)$]  $\dim R'/I(G^{S}) = d$. 
		\item[$(ii)$]  $h_{s}(R'/I(G^{S})) = h_{s}(R/I(G)) + (-1)^{s - 1}$. 
		%where  $h_{R'/I(G^{S})}(\lambda) = h'_0 + h'_1\lambda + h'_2\lambda^2 + \cdots + h'_s\lambda^s$.
	\end{enumerate}
%	$h'_{s} = h_{s} + (-1)^{s - 1}$, where  
%	$h_{R'/I(G^{S})}(\lambda) = h'_0 + h'_1\lambda + h'_2\lambda^2 + \cdots + h'_s\lambda^s$. 
%	\begin{enumerate}
%		\item[$(i)$]  
%	\end{enumerate}
\end{enumerate}
\end{Lemma}
\begin{proof}
\begin{enumerate}
	\item Since $G$ is an induced subgraph of $G^{S}$, we have $\im(G^{S}) \geq \im(G)$. 
	Put $t = \im(G)$. 
	Assume that $t < \im(G^{S})$.  
	Then there is an induced matching of $G^{S}$ consists of $t + 1$ edges 
	$\{x_{i_{1}}, x_{j_{1}}\}$, \ldots, $\{x_{i_{t}}, x_{j_{t}}\}$ and $\{x_{i_{t + 1}}, x_{j_{t + 1}}\}$, 
	where $1 \leq i_{1}, \ldots, i_{t + 1}, j_{1}, \ldots, j_{t + 1} \leq n + 1$. 
	It cannot be an induced matching of $G$ since $t = \im(G)$. 
	Hence we may assume $x_{j_{t + 1}} = x_{n + 1}$. 
	Then $\{ x_{i_{1}}, x_{j_{1}} \} \in E(G)$. 
	By definition of induced matching, one has $\{x_{i_{1}}, x_{n + 1}\} \not\in E(G^{S})$ and 
	$\{x_{j_{1}}, x_{n + 1}\} \not\in E(G^{S})$. 
	It follows that $x_{i_{1}}, x_{j_{1}} \in S$, but this is a contradiction since $S$ is an 
	independent set of $G$. 
	Therefore one has $\im(G^{S}) = \im(G)$. 
	\item Since $G$ is an induced subgraph of $G^{S}$, one has $\reg(I(G)) \leq \reg(I(G^{S}))$ 
	by \cite[Lemma 3.1(1)]{Ha}.  
	Note that $I(G^{S}) + (x_{n + 1}) = (x_{n + 1}) + I(G)$ and 
	$I(G^{S}) : (x_{n + 1}) = (x_{i} : x_{i} \not\in S)$. 
	Hence $\reg (I(G)^{S} + (x_{n + 1})) = \reg (I(G))$ and 
	$\reg (I(G^{S}) : (x_{n + 1})) = 1$ by Lemma \ref{LemB} (2). 
	We also note that $\reg (I(G)) \geq 2$ since $G$ has no isolated vertex. 
	Thus one has $\reg (I(G^{S})) = \reg (I(G))$ by Lemma \ref{LemC}. 
	Therefore $\reg(R'/I(G^{S})) = \reg(R/I(G))$ by Lemma \ref{LemA}. 
	\item By additivity of Hilbert series on the exact sequence 
	\[
	0 \to R'/I(G^{S}) : (x_{n + 1}) \ (-1) \xrightarrow{\ \cdot x_{n + 1} \ }  R'/I(G^{S}) \to R'/I(G^{S}) + (x_{n + 1}) \to 0, 
	\]
	it follows that $H_{R'/I(G^{S})}(\lambda) = H_{R'/I(G^{S}) + (x_{n + 1})}(\lambda) + \lambda \cdot H_{R'/I(G^{S}) : (x_{n + 1})}(\lambda)$. 
	As seen before, $I(G^{S}) + (x_{n + 1}) = (x_{n + 1}) + I(G)$ and 
	$I(G^{S}) : (x_{n + 1}) = (x_{i} : x_{i} \not\in S)$. 
	Hence we have $H_{R'/I(G^{S}) + (x_{n + 1})}(\lambda) = H_{R/I(G)}(\lambda)$ and 
	$H_{R'/I(G^{S}) : (x_{n + 1})}(\lambda) = 1/(1 - \lambda)^{|S| + 1}$. 
	Therefore we have the desired conclusion. \qed
\end{enumerate}
\end{proof}

\begin{Example}\normalfont
Let $G = C_{5}$ be the pentagon graph which appears in Example \ref{C5}. 
By Lemma \ref{S-suspension}, we have $\im(G^{\emptyset}) = \im(G) = 1$, $\reg (R'/I(G^{\emptyset})) = \reg (R/I(G)) = 2$ 
and 
\[
\displaystyle H_{R'/I(G^{\emptyset})}(\lambda) = H_{R/I(G)}(\lambda) + \frac{\lambda(1 - \lambda)}{(1 - \lambda)^{2}} = \frac{1 + 3\lambda + \lambda^{2}}{(1 - \lambda)^{2}} + \frac{\lambda(1 - \lambda)}{(1 - \lambda)^{2}} = \frac{1 + 4\lambda}{(1 - \lambda)^{2} }.
\] 
\end{Example}

%Let $G$ be a finite simple graph on $V(G) = \{ x_{1}, \ldots, x_{n} \}$. 
%For $x \in V(G)$, let $N_{G}(x)$ denote the neighborhood of $x$. 
%For a subset $A \subset V(G)$, let $N_{G}(A) = \bigcup_{x \in A} N_{G}(x) \setminus A$. 

Let $\{ x_{i}, x_{j} \} \in E(G)$ be an edge of $G$ and let $S \subset V(G)$ be an independent set of $G$ 
such that $\{x, x_{i}\} \not\in E(G)$ and $\{x, x_{j}\} \not\in E(G)$ for all $x \in S$. 
The next lemma mentions the induced matching number, regularity and Hilbert series of 
the $\{x_{i}, x_{j}\}, S$-suspension of $G$. \\

\begin{Lemma}
\label{e,S-suspension}
Use the notation as above. 
Suppose that, for all $x \in V(G) \setminus \left( S \cup \{x_{i}\} \cup \{x_{j}\} \right)$, there exists $x_{k} \in S \cup \{x_{i}\} \cup \{x_{j}\}$ 
such that $\{x, x_{k}\} \in E(G)$. 
%Let $G' = G^{\{x_{i}, x_{j}\}, S}$.  
Then one has
	\item[$(1)$]  $\im(G^{\{x_{i}, x_{j}\}, S}) = \im(G)$. 
	\item[$(2)$]  $\reg(R'/I(G^{\{x_{i}, x_{j}\}, S})) = \reg(R/I(G))$ if $\reg(R/I(G)) \geq 2$. 
	\item[$(3)$]  The Hilbert series of $R'/I(G^{\{x_{i}, x_{j}\}, S})$ is 
	\begin{eqnarray*} H_{R'/I(G^{\{x_{i}, x_{j}\}, S})}(\lambda) &=& H_{R/I(G)}(\lambda) + \frac{\lambda(1 + \lambda)}{(1 - \lambda)^{|S| + 2}} \\
	&=& \begin{cases} 
		\displaystyle \frac{h_{R/I(G)}(\lambda) + \lambda(1 + \lambda)(1 - \lambda)^{d - |S| - 2} }{(1 - \lambda)^{d}} & \text{$(|S| \leq d - 2)$}, \\
		\displaystyle \frac{(1 - \lambda)^{|S| + 2 - d}h_{R/I(G)}(\lambda) + \lambda(1 + \lambda)}{(1 - \lambda)^{|S| + 2}} & \text{$(d - 1 \leq |S| \leq d)$}. 
		\end{cases} 
	\end{eqnarray*}
	In particular, if $|S| = d - s$ and $s \geq 2$, we have 
	\begin{enumerate}
		\item[$(i)$]  $\dim R'/I(G^{\{x_{i}, x_{j}\}, S}) = d$. 
		\item[$(ii)$]  $h_{s}(R'/I(G^{\{x_{i}, x_{j}\}, S})) = h_{s}(R/I(G)) + (-1)^{s - 2}$. 
%		, where  
%	$h_{R'/I(G')}(\lambda) = h'_0 + h'_1\lambda + h'_2\lambda^2 + \cdots + h'_s\lambda^s$.
	\end{enumerate}
\end{Lemma}
\begin{proof}
\begin{enumerate}
	\item Since $G$ is an induced subgraph of $G^{\{x_{i}, x_{j}\}, S}$, we have $\im(G^{\{x_{i}, x_{j}\}, S}) \geq \im(G)$. 
	Put $t = \im(G)$. 
	Assume that $t < \im(G^{\{x_{i}, x_{j}\}, S})$.  
	Then there is an induced matching of $G^{\{x_{i}, x_{j}\}, S}$ consists of $t + 1$ edges 
	$\{x_{i_{1}}, x_{j_{1}}\}$, \ldots, $\{x_{i_{t}}, x_{j_{t}}\}$ and $\{x_{i_{t + 1}}, x_{j_{t + 1}}\}$, 
	where $1 \leq i_{1}, \ldots, i_{t + 1}, j_{1}, \ldots, j_{t + 1} \leq n + 1$. 
	It cannot be an induced matching of $G$ since $t = \im(G)$. 
	Hence we may assume $x_{j_{t + 1}} = x_{n + 1}$. 
	By the definition of $G^{\{x_{i}, x_{j}\}, S}$, $x_{i_{k}}, x_{j_{k}} \in S \cup \{x_{i}\} \cup \{x_{j}\}$ for all $1 \leq k \leq t$. 
	Since $S$ is an independent set, $x_{i_{k}} \in \{x_{i}\} \cup \{x_{j}\}$ or $x_{j_{k}} \in \{x_{i}\} \cup \{x_{j}\}$ 
	for all $1 \leq k \leq t$. Thus we have $t = 1$ and we may assume that $x_{i_{1}} = x_{i}$ and $x_{j_{1}} = x_{j}$. 
	Since $\{x_{i_{2}}, x_{n + 1}\} \in E(G^{\{x_{i}, x_{j}\}, S})$, $x_{i_{2}} \not\in S \cup \{x_{i}\} \cup \{x_{j}\}$. 
	Hence there exists $x_{k} \in S \cup \{x_{i}\} \cup \{x_{j}\}$ such that $\{ x_{i_{2}}, x_{k} \} \in E(G)$.  
	Then $x_{k} \in S$ since $x_{k} \neq x_{i}$ and $x_{k} \neq x_{j}$. 
	Thus $\{x_{k}, x_{i}\} \not\in E(G)$ and $\{x_{k}, x_{j}\} \not\in E(G)$ from the assumption. 
	This facts says that $\left\{ \{x_{i}, x_{j}\}, \{x_{i_{2}}, x_{k}\} \right\}$ is an induced matching of $G$, 
	contradicts that $\im(G) = t = 1$. 
	Therefore one has $\im(G^{\{x_{i}, x_{j}\}, S}) = \im(G)$. 
	\item Since $G$ is an induced subgraph of $G^{\{x_{i}, x_{j}\}, S}$, 
	one has $\reg(I(G)) \leq \reg(I(G^{\{x_{i}, x_{j}\}, S}))$ by \cite[Lemma 3.1(1)]{Ha}.  
	We note that $I(G^{\{x_{i}, x_{j}\}, S}) + (x_{n + 1}) = (x_{n + 1}) + I(G)$ and 
	$I(G^{\{x_{i}, x_{j}\}, S}) : (x_{n + 1}) = (x_{i}x_{j}) + \left(x_{k} : x_{k} \not\in S \cup \{x_{i}\} \cup \{x_{j}\} \right)$. 
	Hence it follows that $\reg (I(G^{\{x_{i}, x_{j}\}, S}) + (x_{n + 1})) = \reg (I(G))$ and 
	$\reg (I(G^{\{x_{i}, x_{j}\}, S}) : (x_{n + 1})) = 2$ by Lemma \ref{LemB} (2). 
	By virtue of Lemma \ref{LemA}, $\reg (I(G)) \geq 3$ if $\reg (R/I(G)) \geq 2$.  
	Thus one has $\reg (I(G^{\{x_{i}, x_{j}\}, S})) = \reg (I(G))$ by Lemma \ref{LemC}. 
	Therefore $\reg(R'/I(G^{\{x_{i}, x_{j}\}, S})) = \reg(R/I(G))$ by Lemma \ref{LemA}. 
	\item By additivity of Hilbert series on the exact sequence 
	\[
	0 \to R'/I(G^{\{x_{i}, x_{j}\}, S}) : (x_{n + 1}) \ (-1) \xrightarrow{\ \cdot x_{n + 1} \ }  R'/I(G^{\{x_{i}, x_{j}\}, S}) \to R'/I(G^{\{x_{i}, x_{j}\}, S}) + (x_{n + 1}) \to 0, 
	\]
	one has $H_{R'/I(G^{\{x_{i}, x_{j}\}, S})}(\lambda) = H_{R'/I(G^{\{x_{i}, x_{j}\}, S}) + (x_{n + 1})}(\lambda) + \lambda \cdot H_{R'/I(G^{\{x_{i}, x_{j}\}, S}) : (x_{n + 1})}(\lambda)$. 
	Since $I(G^{\{x_{i}, x_{j}\}, S}) + (x_{n + 1}) = (x_{n + 1}) + I(G)$ and 
	$I(G^{\{x_{i}, x_{j}\}, S}) : (x_{n + 1}) = (x_{i}x_{j}) + \left(x_{k} : x_{k} \not\in S \cup \{x_{i}\} \cup \{x_{j}\} \right)$, 
	it follows that  $H_{R'/I(G^{\{x_{i}, x_{j}\}, S}) + (x_{n + 1})}(\lambda) = H_{R/I(G)}(\lambda)$ and 
	$H_{R'/I(G^{\{x_{i}, x_{j}\}, S}) : (x_{n + 1})}(\lambda) = (1 + \lambda)/(1 - \lambda)^{|S| + 2}$. 
	Therefore we have the desired conclusion. \qed
\end{enumerate}
\end{proof}

\begin{Remark}\normalfont
In general, both of the $S$-suspension and $\{x_{i}, x_{j}\}, S$-suspension do not preserve the matching number. 
In fact, let $G = C_{5}, G^{\emptyset}$ and $G^{\{x_{1}, x_{2}\}, \emptyset}$ be the graphs which appears in Example  \ref{C5}. 
Then $\m(G) = 2$ but $\m(G^{\emptyset}) = \m(G^{\{x_{1}, x_{2}\}, \emptyset}) = 3$. 
\end{Remark}

By combining Lemma \ref{S-suspension} and \ref{e,S-suspension}, one has

\begin{Lemma}
\label{h_s}
Let $G$ be a finite simple graph on the vertex set $V(G) = \{x_{1}, \ldots, x_{n}\}$
such that $G$ has no isolated vertices. 
Assume that $\reg (R/I(G)) \geq 2$,
and $h_{R/I(G)}(\lambda) = h_{0}(R/I) + h_{1}(R/I)\lambda + h_{2}(R/I)\lambda^2 + \cdots + h_{s}(R/I)\lambda^s$ with $h_{s}(R/I) \neq 0$.
Then one has 
\begin{enumerate}
	\item[$(1)$] There exists a connected graph $G_{1}$ such that 
	\begin{itemize}
		\item $\im(G_{1}) = \im(G)$. 
		\item $\reg (K[V(G_{1})]/I(G_{1})) = \reg (R/I(G))$. 
		\item $\deg h_{K[V(G_{1})]/I(G_{1})}(\lambda) = \deg h_{R/I(G)}(\lambda) + 1$. 
	\end{itemize}
	\item[$(2)$] If $s \geq 2$, then there exists a connected graph $G_{2}$ such that 
	\begin{itemize}
		\item $\im(G_{2}) = \im(G)$. 
		\item $\reg (K[V(G_{2})]/I(G_{2})) = \reg (R/I(G))$. 
		\item $\deg h_{K[V(G_{2})]/I(G_{2})}(\lambda) = \deg h_{R/I(G)}(\lambda) - 1$. 
	\end{itemize}
\end{enumerate}
\end{Lemma}
\begin{proof} Write $d = \dim R/I(G)$. 
\begin{enumerate}
	\item Take an independent set $S \subset V(G)$ of $G$ with $|S| = d$. 
	Then Lemma \ref{S-suspension} says that $G_{1} = G^{S}$ is the desired graph. 
	\item We distinguish the following four cases. 
	\begin{itemize}
			\item {\bf The case $h_{s} > 0$ and $s$ is even : } Take an independent set $S \subset V(G)$ of $G$ with $|S| = d - s$. 
			By virtue of Lemma \ref{S-suspension}, it follows that 
			\begin{itemize}
				\item $\im(G^{S}) = \im(G)$, 
				\item $\reg (K[V(G^{S})]/I(G^{S})) = \reg (R/I(G))$, 
				\item $\dim K[V(G^{S})]/I(G^{S}) = d$, 
				\item $h_{s}(K[V(G^{S})]/I(G^{S})) = h_{s}(R/I(G)) - 1$, 
			\end{itemize}
			where $K[V(G^{S})] = R \otimes_{K} K[x_{n + 1}]$.   
			By using Lemma \ref{S-suspension} repeatedly, we can construct the desired graph $G_{2}$. 
			\item {\bf The case $h_{s} < 0$ and $s$ is odd : } From the same argument as the case that $h_{s} > 0$ and $s$ is even, 
			we can construct the desired graph $G_{2}$. 
			\item {\bf The case $h_{s} > 0$ and $s$ is odd : } %In this case, note that $s \geq 3$. 
			Since $d - s + 2 \leq d$,  there exists an independent set $T \subset V(G)$ with 
			$|T| = d - s + 2$. 
			Choose two vertices $x_{a}, x_{b} \in T$ and let $S = T \setminus \left(\{x_{a}\} \cup \{x_{b}\}\right)$. 
			Since $|S| = d - s$, by Lemma \ref{S-suspension}, one has 
			\begin{itemize}
				\item $\im(G^{S}) = \im(G)$,
				\item $\reg (K[V(G^{S})]/I(G^{S})) = \reg (R/I(G))$, 
				\item $\dim K[V(G^{S})]/I(G^{S}) = d$, 
				\item $\deg h_{K[V(G^{S})]/I(G^{S})}(\lambda) = s$ and $h_{s}(K[V(G^{S})]/I(G^{S})) = h_{s}(R/I(G)) + 1$.						\end{itemize}
			
			$\ $
				
			 Next, we consider the $\{x_{a}, x_{n + 1}\}, S$-suspension of $G^{S}$. Let $G' = (G^{S})^{\{x_{a}, x_{n + 1}\}, S}$. Note that $\{ x, x_{a} \} \not\in E(G^{S})$ and $\{ x, x_{n + 1} \} \not\in E(G^{S})$ for all $x \in S$. 
			 In addition, for all $x \in V(G^{S}) \setminus \left( S \cup \{x_{a}\} \cup \{x_{n + 1}\} \right)$, $\{x, x_{n + 1}\} \in E(G^{S})$ since $x \not\in S$. Hence, by virtue of Lemma \ref{e,S-suspension}, we have 
			 \begin{itemize}
				\item $\im(G') = \im(G^{S})$,
				\item $\reg (K[V(G')]/I(G')) = \reg (K[V(G^{S})]/I(G^{S}))$, 
				%since $\reg (K[V(G^{S})]/I(G^{S})) = \reg (R/I(G)) \geq 2$, 
				\item $\dim K[V(G')]/I(G') = \dim K[V(G^{S})]/I(G^{S})$, 
				\item $\deg h_{K[V(G')]/I(G')}(\lambda) = \deg h_{K[V(G^{S})]/I(G^{S})}(\lambda) = s$ and \\ $h_{s}(K[V(G')]/I(G')) = h_{s}(K[V(G^{S})]/I(G^{S})) - 1 = h_{s}(R/I(G))$, 
			\end{itemize}
			where $K[V(G')] = K[V(G^{S})] \otimes_{K} K[x_{n + 2}]$. 
			
			$\ $
			
			Finally, we consider the $\{x_{b}, x_{n + 1}\}, S$-suspension of $G'$. 
			Let $G'' = (G')^{\{x_{b}, x_{n + 1}\}, S}$. Note that $\{ x, x_{b} \} \not\in E(G')$ and $\{ x, x_{n + 1} \} \not\in E(G')$ for all $x \in S$.
			In addition, for all $x \in V(G') \setminus \left( S \cup \{x_{b}\} \cup \{x_{n + 1}\} \right)$, $\{x, x_{n + 1}\} \in E(G')$ since $x \not\in S$. Hence, by using Lemma \ref{e,S-suspension} again, we have  
			\begin{itemize}
				\item $\im(G'') = \im(G')$,
				\item $\reg (K[V(G'')]/I(G'')) = \reg (K[V(G')]/I(G'))$, 
				%since $\reg (K[V(G^{S})]/I(G^{S})) = \reg (R/I(G)) \geq 2$, 
				\item $\dim K[V(G'')]/I(G'') = \dim K[V(G')]/I(G')$, 
				\item $\deg h_{K[V(G'')]/I(G''))}(\lambda) = \deg h_{K[V(G')]/I(G'))}(\lambda) = s$ and \\ $h_{s}(K[V(G'')]/I(G'')) = h_{s}(K[V(G')]/I(G')) - 1 = h_{s}(K[V(G)]/I(G)) - 1$ if $h_{s}(K[V(G)]/I(G)) > 1$, 
				\item $\deg h_{K[V(G'')]/I(G''))}(\lambda) < \deg h_{K[V(G')]/I(G'))}(\lambda) = s$ \\ if $h_{s}(K[V(G)]/I(G)) = 1$, 
				%\item $\deg h_{K[V(G'')]/I(G'')}(\lambda) = \deg h_{K[V(G')]/I(G')}(\lambda)$ and $h_{s}(K[V(G'')]/I(G'')) = h_{s}(K[V(G')]/I(G')) - 1$, 
			\end{itemize}
			where $K[V(G'')] = K[V(G')] \otimes_{K} K[x_{n + 3}]$. 
			
			$\ $
			
			Thus, it follows that
			\begin{itemize}
				\item $\im(G'') = \im(G)$,
				\item $\reg (K[V(G'')]/I(G'')) = \reg (R/I(G))$, 
				\item $\dim K[V(G'')]/I(G'') = \dim R/I(G)$, 
				\item $\deg h_{K[V(G'')]/I(G'')}(\lambda) = \deg h_{K[V(G)]/I(G)}(\lambda) = s$ and \\ $h_{s}(K[V(G'')]/I(G'')) = h_{s}(R/I(G)) - 1$ if $h_{s}(R/I(G)) > 1$, 
				\item $\deg h_{K[V(G'')]/I(G'')}(\lambda) < \deg h_{K[V(G)]/I(G)}(\lambda) = s$ if $h_{s}(R/I(G)) = 1$. 
			\end{itemize}
		Therefore, by repeating this argument, we can construct the desired graph $G_{2}$. 
		\item {\bf The case $h_{s} < 0$ and $s$ is even : } From the same argument as the case that $h_{s} > 0$ and $s$ is odd, 
			we can construct the desired graph $G_{2}$. \qed
	\end{itemize}
\end{enumerate}
\end{proof}

\begin{Example}\normalfont
Let $G = C_{8}$ be the octagon graph with $V(G) = \{ x_{1}, x_{2}, \ldots, x_{8} \}$ and 
$E(G) = \left\{ \{x_{1}, x_{2}\}, \{x_{2}, x_{3}\}, \ldots, \{x_{7}, x_{8}\}, \{x_{1}, x_{8}\} \right\}$.  Then $\im(G) = 2$. 
By virtue of \cite[Theorem 5.2]{BHT} and \cite[Proof of Proposition 2.10]{HKMT}, one has $\reg (K[V(G)]/I(G)) = 3$ and $\displaystyle H_{K[V(G)]/I(G)}(\lambda) = \frac{1 + 4\lambda + 2\lambda^{2} - 4\lambda^{3} - \lambda^{4}}{(1 - \lambda)^{4}}$. \\

Let $G' = (G^{\emptyset})^{\{x_{1}, x_{9}\}, \emptyset}$ and $G'' = (G')^{\{x_{3}, x_{9}\}, \emptyset}$. 
By using Lemmata \ref{S-suspension} and \ref{e,S-suspension}, we have 
\begin{itemize} 
	\item $\im(G) = \im(G^{\emptyset}) = \im(G') = \im(G'') = 2$, 
	\item $\reg(K[V(G)]/I(G)) = \reg(K[V(G^{\emptyset})]/I(G^{\emptyset})) = \reg(K[V(G')]/I(G'))$ \\  $= \reg(K[V(G'')]/I(G'')) = 3$, 
	\item $\deg h_{K[V(G)]/I(G)}(\lambda) = \deg h_{K[V(G^{\emptyset})]/I(G^{\emptyset})}(\lambda) = \deg h_{K[V(G')]/I(G')}(\lambda) = 4$ and \\ $\deg h_{K[V(G'')]/I(G'')}(\lambda) = 3$. 
\end{itemize} 
Moreover, the Hilbert series of these graphs are as follows. 

\bigskip

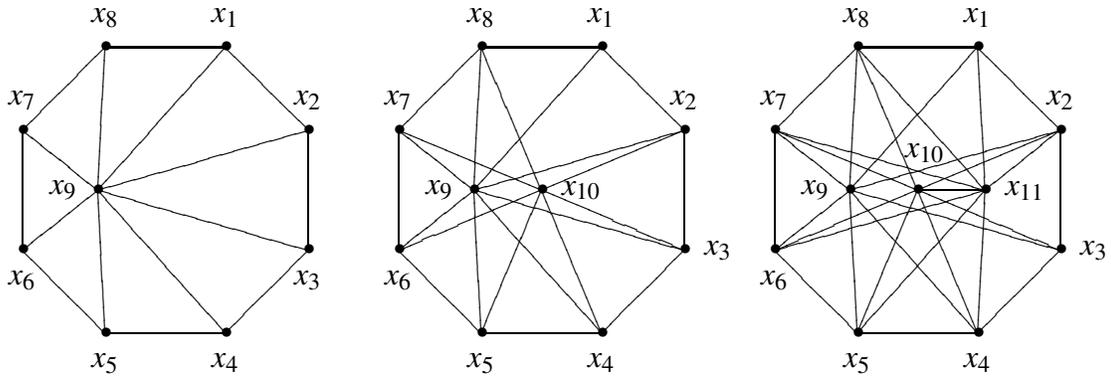
\begin{xy}
	\ar@{} (0,0);(20, 0) *++!D{x_{8}} *\dir<4pt>{*} = "H1";
	\ar@{-} "H1";(36, 0) *++!D{x_{1}} *\dir<4pt>{*} = "A1";
	\ar@{-} "A1";(47, -11) *++!D{x_{2}} *\dir<4pt>{*} = "B1";
	\ar@{-} "B1";(47, -27) *++!U{x_{3}} *\dir<4pt>{*} = "C1";
	\ar@{-} "C1";(36, -38) *++!U{x_{4}} *\dir<4pt>{*} = "D1";
	\ar@{-} "D1";(20, -38) *++!U{x_{5}} *\dir<4pt>{*} = "E1";
	\ar@{-} "E1";(9, -27) *++!U{x_{6}} *\dir<4pt>{*} = "F1";
	\ar@{-} "F1";(9, -11) *++!D{x_{7}} *\dir<4pt>{*} = "G1";
	\ar@{-} "G1";"H1";
	\ar@{-} "A1";(19, -19) *++!R{x_{9}} *\dir<4pt>{*} = "I1";
	\ar@{-} "B1";"I1";
	\ar@{-} "C1";"I1";
	\ar@{-} "D1";"I1";
	\ar@{-} "E1";"I1";
	\ar@{-} "F1";"I1";
	\ar@{-} "G1";"I1";
	\ar@{-} "H1";"I1";
	\ar@{} (0,0);(70, 0) *++!D{x_{8}} *\dir<4pt>{*} = "H2";
	\ar@{-} "H2";(86, 0) *++!D{x_{1}} *\dir<4pt>{*} = "A2";
	\ar@{-} "A2";(97, -11) *++!D{x_{2}} *\dir<4pt>{*} = "B2";
	\ar@{-} "B2";(97, -27) *++!L{x_{3}} *\dir<4pt>{*} = "C2";
	\ar@{-} "C2";(86, -38) *++!U{x_{4}} *\dir<4pt>{*} = "D2";
	\ar@{-} "D2";(70, -38) *++!U{x_{5}} *\dir<4pt>{*} = "E2";
	\ar@{-} "E2";(59, -27) *++!U{x_{6}} *\dir<4pt>{*} = "F2";
	\ar@{-} "F2";(59, -11) *++!D{x_{7}} *\dir<4pt>{*} = "G2";
	\ar@{-} "G2";"H2";
	\ar@{-} "A2";(69, -19) *++!R{x_{9}} *\dir<4pt>{*} = "I2";
	\ar@{} "I2";(78, -19) *++!L{x_{10}} *\dir<4pt>{*} = "J2";
	\ar@{-} "B2";"I2";
	\ar@{-} "C2";"I2";
	\ar@{-} "D2";"I2";
	\ar@{-} "E2";"I2";
	\ar@{-} "F2";"I2";
	\ar@{-} "G2";"I2";
	\ar@{-} "H2";"I2";
	\ar@{-} "B2";"J2";
	\ar@{-} "C2";"J2";
	\ar@{-} "D2";"J2";
	\ar@{-} "E2";"J2";
	\ar@{-} "F2";"J2";
	\ar@{-} "G2";"J2";
	\ar@{-} "H2";"J2";
	\ar@{} (0,0);(120, 0) *++!D{x_{8}} *\dir<4pt>{*} = "H3";
	\ar@{-} "H3";(136, 0) *++!D{x_{1}} *\dir<4pt>{*} = "A3";
	\ar@{-} "A3";(147, -11) *++!D{x_{2}} *\dir<4pt>{*} = "B3";
	\ar@{-} "B3";(147, -27) *++!L{x_{3}} *\dir<4pt>{*} = "C3";
	\ar@{-} "C3";(136, -38) *++!U{x_{4}} *\dir<4pt>{*} = "D3";
	\ar@{-} "D3";(120, -38) *++!U{x_{5}} *\dir<4pt>{*} = "E3";
	\ar@{-} "E3";(109, -27) *++!U{x_{6}} *\dir<4pt>{*} = "F3";
	\ar@{-} "F3";(109, -11) *++!D{x_{7}} *\dir<4pt>{*} = "G3";
	\ar@{-} "G3";"H3";
	\ar@{-} "A3";(119, -19) *++!R{x_{9}} *\dir<4pt>{*} = "I3";
	\ar@{} "I3";(128, -19) *\dir<4pt>{*} = "J3";
	\ar@{-} "J3";(137, -19) *++!L{x_{11}} *\dir<4pt>{*} = "K3";
	\ar@{-} "B3";"I3";
	\ar@{-} "C3";"I3";
	\ar@{-} "D3";"I3";
	\ar@{-} "E3";"I3";
	\ar@{-} "F3";"I3";
	\ar@{-} "G3";"I3";
	\ar@{-} "H3";"I3";
	\ar@{-} "B3";"J3";
	\ar@{-} "C3";"J3";
	\ar@{-} "D3";"J3";
	\ar@{-} "E3";"J3";
	\ar@{-} "F3";"J3";
	\ar@{-} "G3";"J3";
	\ar@{-} "H3";"J3";
	\ar@{-} "A3";"K3";
	\ar@{-} "B3";"K3";
	\ar@{-} "D3";"K3";
	\ar@{-} "E3";"K3";
	\ar@{-} "F3";"K3";
	\ar@{-} "G3";"K3";
	\ar@{-} "H3";"K3";
	\ar@{} (0,0);(129,-14) *{\text{$x_{10}$}};
\end{xy}

\bigskip

\hspace{20mm} $G^{\emptyset} \xrightarrow{\ \ \ \ \ \ \ \{x_{1}, x_{9}\}, \emptyset \text{-suspension}\ \ \ \ \ \ \ } G' \xrightarrow{\ \ \ \ \ \ \ \{x_{3}, x_{9}\}, \emptyset \text{-suspension}\ \ \ \ \ \ \ } G''$

\bigskip

$\displaystyle \frac{1 + 5\lambda - \lambda^{2} - \lambda^{3} - 2\lambda^{4}}{(1 - \lambda)^{4}} \xrightarrow{\ \ \ \ \ \ \ } \frac{1 + 6\lambda - 2\lambda^{2} - 2\lambda^{3} - \lambda^{4}}{(1 - \lambda)^{4}} \xrightarrow{\ \ \ \ \ \ \ } \frac{1 + 7\lambda - 3\lambda^{2} - 3\lambda^{3}}{(1 - \lambda)^{4}}$

\bigskip

\end{Example}

We have finished preparing for the proof of the following lemma. \\

%\begin{Lemma}
%\label{gap-free2}
%Let $r, s \geq 1$ be integers. Then there is a finite connected simple graph $G$ 
%with $im(G) = 1$, $\reg (R/I(G)) = r$ and $\deg h_{R/I(G)} (\lambda) = s$. 
%\end{Lemma}
%\begin{proof}
%Let $r, s \geq 1$ be integers.
%
%First, assume that $1 \leq r \leq s$. By Lemma \ref{gap-free1}, there is a finite connected simple graph $L_{r}$ with $im(L_{r}) = 1$ and $\reg (R/I(L_{r})) = \dim R/I(L_{r}) = r$. 
%Since $\dim R/I(L_{r}) = r$, there is an independent set $S \subset V(L_{r})$ with $|S| = r$. 
%Let $L_{r}^{S}$ be the $S$-suspension of $L_{r}$. By Lemma \ref{S-suspension},  
%$L_{r}^{S}$ is a finite connected simple graph with $im(L_{r}^{S}) = 1$, $\reg (R'/I(L_{r}^{S})) = r$ 
%and $\deg h_{R'/I(L_{r}^{S})} = \deg h_{R/I(L_{r})} + 1$, where $R' = R \otimes_{K} K[x_{n + 1}]$.  
%By repeating the similar operations, we can construct a finite connected simple graph $G$ 
%with $im(G) = 1$, $\reg (R/I(G)) = r$ and $\deg h_{R/I(G)} (\lambda) = s$. 
%
%Next, assume that $1 \leq s < r$. 
%Let $h_{R/I(L_{r})} = h_0 + h_1\lambda + h_2\lambda^2 + \cdots + h_{s'} \lambda^{s'}$. Then $s' \leq r$. 
%By using Lemma \ref{S-suspension} and \ref{G'''} repeatedly, we can construct a finite connected simple graph $G$ with $im(G) = 1$, $\reg (R/I(G)) = r$ and $\deg h_{R/I(G)} (\lambda) = s$.  
%\end{proof}

\begin{Lemma}
\label{disconnected}
Let $G$ be a finite simple graph and $G_{1}, \ldots, G_{\ell}$ the connected components of $G$. 
Assume that $|V(G_{i})| \geq 2$ for all $1 \leq i \leq \ell$.  Then one has 
\begin{enumerate}
	\item[$(1)$] $\im(G) = \sum_{i = 1}^{\ell} \im(G_{i})$. 
	\item[$(2)$] $\reg(R/I(G)) = \sum_{i = 1}^{\ell} \reg(R_{i}/I(G_{i}))$. 
	\item[$(3)$] $\dim R/I(G) = \sum_{i = 1}^{\ell} \dim R_{i}/I(G_{i})$.
	\item[$(4)$] $\deg h_{R/I(G)} (\lambda)= \sum_{i = 1}^{\ell} \deg h_{R_{i}/I(G_{i})} (\lambda)$. 
\end{enumerate}
Here $R_{i} = K\left[ x_{j}^{(i)} : x_{j}^{(i)} \in V(G_{i}) \right]$ for all $1 \leq i \leq \ell$ and $R = R_{1} \otimes_{K} \cdots \otimes_{K} R_{\ell}$. 
\end{Lemma}  

%Now we are in the position to prove Theorem \ref{main}. 

\section{Proof of Theorem \ref{main}}

In this section, we give a proof of Theorem \ref{main}. 

\begin{proof}({\em Proof of Theorem \ref{main}.})
Let $1 \leq a \leq r$ and $s \geq 1$ be integers. 

\begin{itemize}
	\item {\bf The case $a = r = 1$ : } Let $G^{{\rm star} (x_{v})}_{s}$ be the star graph which appears in \cite[Figure 2]{HKMT}. Then \cite[Lemma 1.7]{HKMT} says that $G^{{\rm star} (x_{v})}_{s}$ is the desired graph. 
	\item {\bf The case $a = 1$ and $r \geq 2$ : } By virtue of \cite[Theorem 6]{CKV}, there is a finite connected simple graph $G$ 
with $\im(G) = 1$ and $\reg (K[V(G)]/I(G)) = \dim K[V(G)]/I(G) = r$. Hence we can construct the desired graph $G(1, r, s)$ by using Lemma \ref{h_s} repeatedly. 
	\item {\bf The case $2 \leq a \leq r$ : } By virtue of \cite[Theorem 6]{CKV}, there is a finite connected simple graph $L_{r - a + 1}$ 
	with $\im(L_{r - a + 1}) = 1$ and $\reg (K[V(L_{r - a + 1})]/I(L_{r - a + 1})) = \dim K[V(L_{r - a + 1})]/I(L_{r - a + 1}) = r - a + 1$. 
	Let $G$ be the union of $L_{r - a + 1}$ and $a-1$ disjoint edges $\{ x_{i_{1}}, x_{j_{1}} \}, \ldots, \{ x_{i_{a - 1}}, x_{j_{a - 1}} \}$. 		By virtue of Lemma \ref{disconnected}, we have 
	\begin{itemize}
		\item $\im(G) = 1 + (a - 1) = a$, 
		\item $\reg (K[V(G)]/I(G)) = (r - a + 1) + (a - 1) = r$, 
		\item $\dim K[V(G)]/I(G) = (r - a + 1) + (a - 1) = r$, 
		\item $\deg h_{K[V(G)]/I(G)}(\lambda) = \deg h_{K[V(L_{r - a + 1})]/I(L_{r - a + 1})}(\lambda) + a - 1$. 
	\end{itemize}
	Hence, by using Lemma \ref{h_s} repeatedly, we can construct the desired graph $G(a, r, s)$. \qed 
\end{itemize}

%Let $G_{0}$ be the union of $L_{r - a + 1}$ and $a-1$ disjoint edges $\{ x_{i_{1}}, x_{j_{1}} \}, \ldots, \{ x_{i_{a - 1}}, x_{j_{a - 1}} \}$. By virtue of Lemma \ref{disconnected}, we have 
%\[
%im(G_{0}) = 1 + (a - 1) = a, 
%\]
%\[
%\reg (R/I(G_{0})) = (r - a + 1) + (a - 1) = r, 
%\]
%\[
%\dim R/I(G_{0}) = (r - a + 1) + (a - 1) = r. 
%\]
%By using Lemma \ref{S-suspension} and \ref{G'''} repeatedly, we can construct a finite connected simple graph $G$ with $im(G) = a$, $\reg (R/I(G)) = r$ and $\deg h_{R/I(G)} (\lambda) = s$.  
\end{proof}

\begin{Remark}\normalfont
From Theorem \ref{main} and the inequalities $\im(G) \leq \reg \left(R/I(G)\right) \leq \m(G)$, it is natural to ask the following question : given any integers $a, r, m, s$ with $1 \leq a \leq r \leq m$ 
and $s \geq 1$, is there a simple connected graph $G = G(a, r, m, s)$ with $\im(G) = a$, $\reg (R/I(G)) = r$, $\m(G) = m$ and 
$\deg h_{R/I(G)} = s$ ? 

However, \cite[Theorem 11]{T} says that there is no simple connected graphs $G$ with 
$\im(G) = 1$ and $\reg (R/I(G)) = \m(G) = r$ for all $r \geq 3$. 
\end{Remark}

\bigskip

\noindent
{\bf Acknowledgment.}
The first author was partially supported by JSPS KAKENHI 26220701.
The third author was partially supported by JSPS KAKENHI 17K14165.

\end{document}